\documentclass[11pt,twoside]{article}
\usepackage[T1]{fontenc}
\usepackage{amsmath,amssymb,amsthm,mathtools}
\usepackage{newtxtext,newtxmath}
\usepackage[paperwidth=170mm,paperheight=240mm,
  left=22.5mm,right=22.5mm,top=29mm,bottom=18mm,
  headheight=14pt,headsep=3mm]{geometry}
\usepackage{microtype}
\usepackage{titlesec}
\usepackage{fancyhdr}
\usepackage{etoolbox}
\usepackage{xurl}
\usepackage[colorlinks=true,linkcolor=blue,citecolor=blue,urlcolor=blue]{hyperref}
\newtheoremstyle{paperplain}{10pt}{10pt}{\itshape}{}%
  {\bfseries}{.}{.5em}{}
\newtheoremstyle{paperremark}{10pt}{10pt}{\normalfont}{}%
  {\bfseries}{.}{.5em}{}
\theoremstyle{paperplain}
\newtheorem{theorem}{Theorem}[section]
\newtheorem{lemma}[theorem]{Lemma}
\newtheorem{proposition}[theorem]{Proposition}
\newtheorem{corollary}[theorem]{Corollary}
\newtheorem*{claim}{Claim}
\theoremstyle{paperremark}
\newtheorem*{AIuse}{AI Usage Statement}
\newtheorem{definition}[theorem]{Definition}
\newtheorem{remark}[theorem]{Remark}
\newtheorem{example}[theorem]{Example}

\BeforeBeginEnvironment{lemma}{\begin{samepage}}
\AfterEndEnvironment{lemma}{\end{samepage}}
\BeforeBeginEnvironment{proposition}{\begin{samepage}}
\AfterEndEnvironment{proposition}{\end{samepage}}

\titleformat{\section}[block]{\normalfont\normalsize\bfseries\centering}
  {\thesection.}{.5em}{}
\titlespacing*{\section}{0pt}{25pt plus 3pt minus 2pt}{15pt plus 2pt minus 1pt}
\newcommand{\diver}{\operatorname{div}}
\newcommand{\Lip}{\operatorname{Lip}}
\newcommand{\dd}{\,dV}
\newcommand{\Tr}{\operatorname{Tr}}
\newcommand{\DM}{\mathcal{DM}}
\newcommand{\ds}{\,d\mathcal H^{n-1}}
\newcommand{\norm}[1]{\left\lVert#1\right\rVert}

\newcommand{\paperauthors}{}
\newcommand{\paperrunningauthors}{}
\newcommand{\paperaffiliations}{}
\newcommand{\papertitle}{Local regularity and $p$-harmonic approximation\\
  for weak inverse mean curvature flow}
\newcommand{\papershorttitle}{Regularity and approximation of weak IMCF}
\hypersetup{pdftitle={Local regularity and p-harmonic approximation for weak inverse mean curvature flow},
  pdfauthor={},pdfsubject={Inverse mean curvature flow},bookmarksdepth=1}

\fancypagestyle{firstpage}{\fancyhf{}}
\makeatletter
\newcommand{\makepapertitle}{%
  \thispagestyle{firstpage}%
  \vspace*{12mm}%
  \begin{center}
    {\fontsize{14}{18}\selectfont\bfseries\papertitle\par}%
    \ifx\paperauthors\@empty\else
      \vspace{11pt}{\normalsize\paperauthors\par}%
    \fi
  \end{center}
  \vspace{19pt}}
\renewenvironment{thebibliography}[1]
  {\section*{References}\small
   \list{[\@arabic\c@enumiv]}%
     {\settowidth\labelwidth{[#1]}%
      \leftmargin\labelwidth\advance\leftmargin\labelsep
      \usecounter{enumiv}\setlength{\itemsep}{6pt}
      \setlength{\parsep}{0pt}}
   \sloppy\clubpenalty4000\widowpenalty4000\interlinepenalty10000
   \sfcode`\.=1000\relax}
  {\endlist}
\makeatother

\renewcommand{\paperauthors}{%
  Wei-Chun Chen \qquad Chun-Kai Lien\\[8pt]
  {\small Department of Mathematics, National Taiwan University,\\
  Taipei 10617, Taiwan}\\[4pt]
  {\small
    \href{mailto:b14201028@ntu.edu.tw}{\texttt{b14201028@ntu.edu.tw}}
    \qquad
    \href{mailto:d13221002@ntu.edu.tw}{\texttt{d13221002@ntu.edu.tw}}
  }
}

\renewcommand{\paperrunningauthors}{Wei-Chun Chen and Chun-Kai Lien}

\hypersetup{
  pdfauthor={Wei-Chun Chen, Chun-Kai Lien}
}

\begin{document}
\makepapertitle

\begingroup\small
\noindent\textbf{Abstract.}
We prove local Lipschitz regularity for locally bounded $BV$ solutions
of the calibration equations for inverse mean curvature flow. We also
approximate every prescribed variational weak flow locally by Moser's logarithmic
transforms of positive $p$-harmonic functions. On a fixed smooth domain,
with the trace of a calibrated Lipschitz solution prescribed, both
the uniform error and the absolute difference of total variations are
$O((p-1)^{1/3})$ as $p\downarrow1$. In particular, the approximations
converge strictly in $BV$.

\medskip
\noindent\textbf{Mathematics Subject Classification (2020).}
53E10, 35J92, 49Q20.

\medskip
\noindent\textbf{Keywords.}
Inverse mean curvature flow, $p$-harmonic functions, weak solutions,
calibrations, $BV$ regularity, comparison principles.
\par\endgroup

\section{Introduction}

The use of inverse mean curvature flow (IMCF) in mass inequalities
goes back to Geroch \cite{Geroch} and Jang and Wald \cite{JangWald}.
Moser \cite{Moser} used positive $p$-harmonic functions to construct weak
solutions in the sense of
Huisken and Ilmanen \cite{HIAnnouncement,HI}. Here we start with a
prescribed weak flow $u$ and ask whether the same approximation
recovers it on an interior domain. This local problem, without
assuming properness or an initial set, is posed in
\cite[Question~1.10]{BMRSX}.

Let $(M,g)$ be a smooth Riemannian manifold without boundary, of
dimension $n\geq2$, and let $\Omega\subset M$ be open. For a $BV$
function $u$ with a calibration, identifying an approximation limit
with $u$ requires control of the flux excess against $|Du|$, which may
charge sets of zero volume. Perturbing the exponential flux weights
in opposite directions supplies the nonnegative source terms for the
comparison in Section~\ref{comparisonsection}.

\begin{theorem}\label{bvregularity}
Let $u\in BV_{\mathrm{loc}}(\Omega)\cap L^\infty_{\mathrm{loc}}(\Omega)$.
Suppose that a locally bounded vector field $\nu$, whose divergence
is a locally finite Radon measure, satisfies
\begin{equation}\label{bvcal}
 |\nu|\leq1\quad\text{a.e.},\qquad
 \diver\nu=(\nu,Du)=|Du|
 \quad\text{as measures on }\Omega.
\end{equation}
Then $u$ has a locally Lipschitz representative, and this representative
is a variational weak solution of IMCF. In particular, the singular part of $Du$ vanishes.
\end{theorem}

By the theorem and the converse calibration result
\cite[Lemma~2.3.4]{XuThesis}, the locally bounded $BV_{\mathrm{loc}}$
calibration class coincides, up to equality a.e., with the locally
Lipschitz variational weak solutions. The proof obtains the Lipschitz
representative by identifying $u$ with a locally Lipschitz limit of
logarithmic transforms of positive $p$-harmonic functions.

\begin{remark}
A jump of the evolving regions corresponds to a plateau of the
level-set function, which is compatible with local Lipschitz regularity.
The calibration already rules out the jump part of the distributional derivative by \cite[Lemma~3.3]{CMS}.
The additional conclusions are the vanishing of the Cantor part and
the local gradient bound.
\end{remark}

Here $(\nu,Du)$ is Anzellotti's pairing \cite{Anz}, recalled in
Section~\ref{prelim}. The level-set equation is formally written as
\begin{equation}\label{imcf}
 \diver\!\left(\frac{\nabla u}{|\nabla u|}\right)=|\nabla u|.
\end{equation}

For smooth uniformly convex anisotropies, approximation by
$p$-capacitary potentials was established by Della Pietra, Gavitone,
and Xia \cite{DPGX}. Cabezas-Rivas, Moll, and Solera \cite{CMS}
developed the $BV$ formulation for general norms, including
crystalline anisotropies. Theorem~\ref{bvregularity} starts from
the isotropic calibration equations on an open set.
For proper $BV$ variational flows in dimension three with smooth
ambient metric and $C^2$ initial boundary, regularity is discussed in
\cite[Remark~2.3]{FGP}, where the proof is deferred.

For the approximation statement, write
$\Delta_pv=\diver(|\nabla v|^{p-2}\nabla v)$.
For a positive $p$-harmonic function $v_p$, the logarithmic transform
$u_p=(1-p)\log v_p$ satisfies
\begin{equation}\label{mosertransform}
 \Delta_pu_p=|\nabla u_p|^p.
\end{equation}
We call a weak solution of this equation a $p$-IMCF. The trace of a
calibrated Lipschitz solution determines the limit of the Dirichlet
family on a fixed domain.

\begin{theorem}\label{main}
Let $G\Subset M$ be a smooth domain with nonempty boundary, and let
$u\in W^{1,\infty}(G)\cap C(\overline G)$ admit a calibration on $G$
in the sense of Definition~\ref{caldef}.
For $p>1$, let $v_p\in W^{1,p}(G)\cap C(\overline G)$
be the positive $p$-harmonic Dirichlet solution with boundary
datum $e^{-u/(p-1)}$, and set $u_p=-(p-1)\log v_p$.
There is $C>0$, depending only on $(G,g)$, $\operatorname{osc}_G u$,
and $\|\nabla u\|_{L^\infty(G)}$, such that, for $1<p<2$,
\begin{gather}
 \|u_p-u\|_{C^0(\overline G)}
 \leq C(p-1)^{1/3},\label{mainconvergence}\\
 \left|\int_G|\nabla u_p|\dd-\int_G|\nabla u|\dd\right|
 \leq C(p-1)^{1/3}.\label{mainstrictbvrate}
\end{gather}
Thus the full family converges uniformly to $u$ on $\overline G$
and strictly in $BV(G)$ as $p\downarrow1$.
\end{theorem}

For global constructions and properness see \cite{MRS,BMRSX,XuProper}.
The gap in the Green-kernel iteration of \cite{MRS} and the status of
its conclusions are discussed in
\cite[Sections~3.1 and~3.3]{BMRSX}.
The approximation in \cite[Theorem~2.8]{BPP} starts from a proper flow
issuing from a bounded closed set with $C^{1,1}$ boundary. Strong
gradient convergence in $L^q_{\mathrm{loc}}$, $1\leq q<\infty$, holds
under the further hypotheses of \cite[Theorem~4.9]{BPP}.
Theorem~\ref{main} gives quantitative convergence on the whole fixed
domain when the Dirichlet data come from the prescribed calibrated
solution.

Example~\ref{annulusexample} shows why the trace must come from a
calibrated solution. Restricting the Dirichlet approximations to
interior subdomains gives the local result.

\begin{corollary}\label{localapprox}
For every locally Lipschitz variational weak solution $u$ on $\Omega$
and every domain $D\Subset\Omega$, there are exact $p$-IMCF solutions
on a neighborhood of $\overline D$ converging to $u$ in
$C^0(\overline D)$ as $p\downarrow1$.
\end{corollary}

By Theorem~\ref{bvregularity}, the same conclusion holds for every
locally bounded $BV$ calibration, using its Lipschitz representative.

\begin{AIuse}
The authors formulated the research questions and provided the relevant mathematical background and literature to ChatGPT 6.1 Sol. Through iterative interactions guided by the authors, the AI contributed key ideas and proof strategies for the comparison estimates in Section \ref{comparisonsection}. The authors independently verified and substantially revised the arguments and take full responsibility for the correctness and content of the manuscript.
\end{AIuse}

\section{Weak solutions and Dirichlet approximation}\label{prelim}\label{dirichletsection}

For $K\Subset\Omega$ and locally Lipschitz $u,w$, the functional in
\cite[eq.~(1.7)]{HI} is
\[
 J_u(w;K)=\int_K\bigl(|\nabla w|+w|\nabla u|\bigr)\dd.
\]

\begin{definition}\label{weakdef}
A function $u\in\Lip_{\mathrm{loc}}(\Omega)$ is a weak solution of
\eqref{imcf} if
$J_u(u;K)\leq J_u(w;K)$ whenever
$w\in\Lip_{\mathrm{loc}}(\Omega)$ and
$\{u\ne w\}\Subset K\Subset\Omega$.
\end{definition}

\begin{definition}\label{caldef}
A calibration of a locally Lipschitz function $u$ on $G$ is a
measurable vector field $\nu$ satisfying
\begin{equation}\label{cal}
 |\nu|\leq1,\qquad \nu\cdot\nabla u=|\nabla u|\quad\text{a.e.},
 \qquad \diver \nu=|\nabla u|
\end{equation}
in distributions.
\end{definition}

The calibration formulation appears in \cite[eq.~(3.15)]{HI}.
For $BV$ boundary value problems and comparison, see
\cite[Definition~3.1]{MSBounded} and \cite{MSExterior}.
Moser \cite[Definition~1 and Section~5]{MoserStreamlines} proves the
variational property for locally bounded calibrated
$W^{1,1}_{\mathrm{loc}}$ solutions. For the locally Lipschitz
formulation, see \cite[Definition~2.11]{XuObstacle}.

A locally Lipschitz function with a calibration is a weak solution.
Indeed, if $w-u$ is compactly supported, then \eqref{cal} gives
\begin{align*}
 J_u(w;K)-J_u(u;K)
 &\geq\int_K \nu\cdot\nabla(w-u)+(w-u)|\nabla u|\dd=0.
\end{align*}
Conversely, \cite[Lemma~2.3.4]{XuThesis} gives a calibration for every
locally Lipschitz variational weak solution. We shall use this
equivalence on smooth domains $G\Subset M$.

Denote by $\DM^\infty_{\mathrm{loc}}(G)$ the locally bounded vector fields whose
divergence is a locally finite Radon measure. If
$Z\in\DM^\infty_{\mathrm{loc}}(G)$ and
$v\in BV_{\mathrm{loc}}(G)\cap L^\infty_{\mathrm{loc}}(G)$, its
symmetric pairing with $Dv$ is
\[
 \langle(Z,Dv),\phi\rangle
 =-\int_Gv^*\phi\,d(\diver Z)-\int_GvZ\cdot\nabla\phi\dd.
\]
Here $v^*$ denotes the precise representative; on the jump set it is
the arithmetic mean of the two traces. We use the term $BV$
calibration for \eqref{bvcal}. In the locally Lipschitz class this is
Definition~\ref{caldef}. Recall the standard decomposition of the distributional
derivative of a BV function,
\[
Du=\nabla u\,dV+D^j u+D^c u,
\]
where $D^j u$ and $D^c u$ denote the jump and Cantor
parts, respectively.

The local isotropic argument of
\cite[Lemma~3.3]{CMS} shows that \eqref{bvcal} implies $D^ju=0$ for
every locally bounded $BV$ function. This fact allows us to use the
diffuse-$BV$ chain rule below.

The Green formula, the normal-trace bounds, and the product and chain
rules are recalled in \cite[Sections~2.6--2.8 and Lemma~3.4]{CMS}.
For bounded $BV$ functions $v,a$ with no jump part and
$Z\in\DM^\infty$, we have
\[
 \diver(vZ)=v^*\diver Z+(Z,Dv),\qquad
 (aZ,Dv)=a^*(Z,Dv).
\]
Products with measures are always taken using precise representatives.
For the global Green formula on a smooth domain \(G\), we assume \(Z\in\mathcal{DM}^{\infty}(G)\) and \(v\in BV(G)\cap L^\infty(G)\). Let $n_G$ denote the outward unit normal to $\partial G$ with respect to $g$. The weak normal trace obeys
$\|[Z,n_G]\|_\infty\leq\|Z\|_\infty$, and
\[
 \int_Gv^*\,d(\diver Z)+\int_Gd(Z,Dv)
 =\int_{\partial G}\Tr v\,[Z,n_G]\ds.
\]
Moreover, $[vZ,n_G]=\Tr v\,[Z,n_G]$. On a Riemannian manifold,
these identities follow in local coordinates by including the smooth
volume density in the divergence. The total variation
is understood with respect to $g$.

Let $G\Subset M$ be smooth with nonempty boundary and take
$u\in W^{1,\infty}(G)\cap C(\overline G)$. For $p>1$ the Dirichlet
problem with datum $e^{-u/(p-1)}$ has a unique positive $p$-harmonic
solution $v_p\in W^{1,p}(G)\cap C(\overline G)$ \cite{HKM}.
For $m=\min_{\overline G}u$ and $L=\max_{\overline G}u$, comparison gives
\[
 e^{-L/(p-1)}\leq v_p\leq e^{-m/(p-1)}.
\]
By Moser's logarithmic change of variable \cite[Section~2]{Moser},
the function $u_p=-(p-1)\log v_p$ solves
\begin{equation}\label{pde}
 \begin{cases}
 \Delta_pu_p=|\nabla u_p|^p&\text{in }G,\\
 u_p=u&\text{on }\partial G,
 \end{cases}
\end{equation}
and satisfies
\[
 m\leq u_p\leq L.
\]
Interior $C^{1,\alpha}$ regularity
\cite{Uraltseva,DiBenedetto,Lewis,Tolksdorf}
gives $u_p\in\Lip_{\mathrm{loc}}(G)$.

For $w\in W^{1,p}(G)\cap C(\overline G)$, the change of variable
is reversible:
\[
 \begin{aligned}
 \Delta_pw=|\nabla w|^p
 &\quad\Longleftrightarrow\quad
 \diver(e^{-w}|\nabla w|^{p-2}\nabla w)=0\\
 &\quad\Longleftrightarrow\quad
 \Delta_p(e^{-w/(p-1)})=0.
 \end{aligned}
\]
In the first equivalence the tests $e^{-w}\varphi$ and $e^w\varphi$
admit compactly supported smooth approximations converging in
$W^{1,p}$ and uniformly, as needed for the $L^1$ source term.
The transformation transfers uniqueness and comparison to \eqref{pde}.
In particular, $u\leq v$ on $\partial G$ implies $u_p^u\leq u_p^v$
in $G$, since both changes of variable reverse order. Translating
the datum by $c$ multiplies $v_p$ by $e^{-c/(p-1)}$ and translates
$u_p$ by $c$.

\section{A comparison estimate}\label{comparisonsection}

If $\nu$ calibrates $u$ and $Z=|\nabla U|^{p-2}\nabla U$ is the flux of
a $p$-IMCF $U$, then $e^{-u}\nu$ and $e^{-U}Z$ are divergence-free.
Their difference gives no source term in the Green identity. For
$0<\theta\leq1/2$, put $a=1-\theta$ and $b=1+\theta$. Then
\[
 \diver(e^{-au}\nu-e^{-bU}Z)
 =\theta e^{-au}\,d|Du|+\theta e^{-bU}|\nabla U|^p\,dV.
\]
The paired tests use these measures to control the truncations.

Let $G\Subset M$ be a smooth domain with nonempty boundary. Suppose
that $u\in BV(G)\cap L^\infty(G)$ has a calibration $\nu$, and that
$U\in W^{1,p}(G)\cap\Lip_{\mathrm{loc}}(G)$ is a bounded $p$-IMCF.
We assume $1\leq u,U\leq L$ and use the notation
\[
 \mu=|Du|,\qquad \varepsilon=p-1,\qquad s=|\nabla U|.
\]
Assume that $Z$ is continuous in $G$, as it is for the Dirichlet
solutions used below.

We combine the rescaling of \cite[Theorem~2.2]{HI} with the paired
tests of \cite[Theorem~3.8]{MSBounded}, keeping the boundary terms.
See \cite{CMS} for anisotropic analogues.

Let $\eta$ be a finite nonnegative measure such that
\begin{equation}\label{comparisondefectbounds}
 \begin{aligned}
 (\nu-Z,D(au-bU))&\geq-\eta,\\
 (Z-\nu,D(aU-bu))&\geq-\eta.
 \end{aligned}
\end{equation}
For $k\geq0$, write
\[
 w_\theta=|u-U|-\theta(u+U),\qquad
 E_k=\{w_\theta>k\},\qquad q=\frac{2n}{n-1}.
\]
\begin{proposition}\label{commontruncation}
If $(w_\theta-k)_+$ has compact support in $G$, then
\begin{equation}\label{commontruncationbound}
 \|(w_\theta-k)_+\|_{L^q(G)}^2
 \leq\frac{C}{\theta}\bigl(\eta(E_k)+\varepsilon|E_k|\bigr),
\end{equation}
where $C$ depends only on $(G,g)$ and $L$. If $Z\in L^\infty(G)$,
the same estimate holds without the compact-support assumption, with
the additional term
$\int_{\partial G}\Tr((w_\theta-k)_+)\ds$ on the right inside the
parentheses; in this case $C$ may also depend on $\|Z\|_\infty$.
\end{proposition}
\begin{proof}
We may assume $0\leq k<L$. Put
\[
 E_k^+=\{au-bU>k\},\qquad E_k^-=\{aU-bu>k\}.
\]
On the first set we use the tests
\[
 \chi^+=(au-bU-k)_+,\qquad
 \psi^+=(e^{-bU}-e^{k-au})_+.
\]
Since $D^ju=0$ by \cite[Lemma~3.3]{CMS}, the diffuse-$BV$ chain and
product rules give on $E_k^+$
\[
 D\chi^+=aDu-bDU,\qquad
 D\psi^+=ae^{k-au}Du-be^{-bU}DU.
\]
The calibration and $p$-IMCF equations also give
\[
 \begin{aligned}
 \diver(\nu-Z)&=\mu-s^p\,dV,\\
 \diver(e^{k-au}\nu-e^{-bU}Z)
 &=\theta e^{k-au}\mu+\theta e^{-bU}s^p\,dV.
 \end{aligned}
\]
Test the first flux with $\psi^+$ and the second with $\chi^+$. Let $\mathcal B^+$ denote the sum of the boundary terms arising from the two Green identities.
On \(E_k^+\), the chain rule and the calibration identities give
\[
\begin{aligned}
&(\nu-Z,D\psi^+)
+(e^{k-au}\nu-e^{-bU}Z,D\chi^+)\\
&\quad=(e^{k-au}+e^{-bU})
(\nu-Z,D(au-bU)) -a\psi^+\,d\mu+b\psi^+s^p\,dV.
\end{aligned}
\]
Combining this with the divergence identities above and using \(1-a=b-1=\theta\), the two Green identities yield
\begin{equation}\label{scaledfluxidentity}
 \begin{split}
 \mathcal B^+={}&\int_{E_k^+}(e^{k-au}+e^{-bU})
                   \,d(\nu-Z,D(au-bU))\\
 &+\theta\int_G(\psi^++\chi^+e^{k-au})\,d\mu\\
 &+\theta\int_G(\psi^++\chi^+e^{-bU})s^p\dd.
 \end{split}
\end{equation}
Here $\mathcal B^+=0$ for compactly supported tests. In the boundary case,
\begin{equation}\label{boundaryfluxidentity}
 \begin{split}
 \mathcal B^+={}&\int_{\partial G}\Tr\psi^+\,[\nu-Z,n_G]\ds\\
 &+\int_{\partial G}\Tr\chi^+\,
       [e^{k-au}\nu-e^{-bU}Z,n_G]\ds.
 \end{split}
\end{equation}
The global bound on $Z$ controls both normal traces in
\eqref{boundaryfluxidentity}, with the exponential factors bounded
in terms of $L$.

On $E_k^-$ we use
\[
 \chi^-=(aU-bu-k)_+,\qquad
 \psi^-=(e^{-bu}-e^{k-aU})_+,
\]
testing $Z-\nu$ with $\psi^-$ and $e^{k-aU}Z-e^{-bu}\nu$ with $\chi^-$.
The first divergence is $s^p\,dV-\mu$, and
\[
 \diver(e^{k-aU}Z-e^{-bu}\nu)
 =\theta e^{k-aU}s^p\,dV+\theta e^{-bu}\mu,
\]
so the Green identities give the counterpart of \eqref{scaledfluxidentity}
with $u,U$, $\nu,Z$, and $d\mu,s^p\,dV$ interchanged. Its boundary
term in the boundary case is
\[
 \begin{split}
 \mathcal B^-={}&\int_{\partial G}\Tr\psi^-\,[Z-\nu,n_G]\ds\\
 &+\int_{\partial G}\Tr\chi^-\,
       [e^{k-aU}Z-e^{-bu}\nu,n_G]\ds.
 \end{split}
\]
When $(w_\theta-k)_+$ has compact support, both pairs of tests have
compact support. In that case we apply the identities on a smooth
interior domain and set $\mathcal B^+=\mathcal B^-=0$.

Use the first inequality in \eqref{comparisondefectbounds} on $E_k^+$
and the second on $E_k^-$. The bounds $1\leq u,U\leq L$ and
$0\leq k<L$ keep the exponential weights bounded above and away from
zero in terms of $L$. The nonnegative source terms then give
\[
 \theta\int_G\chi^\pm\,(d\mu+s^p\,dV)
 \leq C\bigl(\eta(E_k^\pm)+|\mathcal B^\pm|\bigr).
\]
In the boundary case, the normal-trace bounds and
$\Tr\psi^\pm\leq C(L)\Tr\chi^\pm$ give
\[
 |\mathcal B^\pm|\leq C(L,\|Z\|_\infty)
                    \int_{\partial G}\Tr\chi^\pm\ds.
\]
The chain rule bounds the variations of $\chi^+$ and $\chi^-$ on their
positivity sets by $a\,d\mu+bs\,dV$ and $b\,d\mu+as\,dV$, respectively.
Using $s\leq s^p+\varepsilon$, $\chi^\pm\leq L$, and $a,b\leq3/2$, we obtain
\[
 \begin{aligned}
 |D((\chi^\pm)^2)|(G)
 &\leq 3\int_G\chi^\pm\,(d\mu+s\,dV)\\
 &\leq\frac C\theta\Bigl(\eta(E_k^\pm)
       +\varepsilon|E_k^\pm|+|\mathcal B^\pm|\Bigr).
 \end{aligned}
\]

Since $E_k=E_k^+\sqcup E_k^-$, we have
$\chi^++\chi^-=(w_\theta-k)_+$ and
$(\chi^+)^2+(\chi^-)^2=(w_\theta-k)_+^2$.
Adding the variation bounds controls the single squared truncation
$v=(w_\theta-k)_+^2$, with the measure terms summed over $E_k$.
In the boundary case the traces add to $\Tr((w_\theta-k)_+)$.
The endpoint $BV$ Sobolev inequality
\[
 \|v\|_{L^{n/(n-1)}(G)}
 \leq C_G\left(|Dv|(G)+\int_{\partial G}|\Tr v|\ds\right)
\]
then gives the required bound for $\|(w_\theta-k)_+\|_{L^q(G)}^2$.
Here the endpoint inequality of \cite{Gagliardo} is applied by
approximation to the zero extension. The boundary contribution follows
from the $BV$ trace and zero-extension formula
\cite[Theorem~3.87]{AFP}.
The trace term vanishes for compactly supported truncations; otherwise
$\Tr((w_\theta-k)_+^2)\leq L\Tr((w_\theta-k)_+)$.
This proves both versions of \eqref{commontruncationbound}.
\end{proof}

\begin{lemma}\label{measuredefect}
For $1<p<2$, the lower bounds in \eqref{comparisondefectbounds} hold with
\begin{equation}\label{bvdefect}
 \eta=C\varepsilon\,dV+C(|Z|-1)_+\,d\mu,
\end{equation}
provided this measure is finite. If $u\in W^{1,\infty}(G)$, one may
instead take
\begin{equation}\label{lipschitzdefect}
 \eta=C_{\|\nabla u\|_\infty}\varepsilon\,dV.
\end{equation}
The constants are independent of $p$ and $\theta$ in the stated ranges.
\end{lemma}
\begin{proof}
Calibration and the pairing bound give
\begin{align*}
 (\nu-Z,D(au-bU))
 &=a\mu+bs^p\,dV-b\nu\cdot\nabla U\,dV-a(Z,Du)\\
 &\geq b(s^p-s)\,dV+a(1-|Z|)\,d\mu.
\end{align*}
The bound $(Z,Du)\leq|Z|\mu$ uses the continuous representative of $Z$
on the Cantor part of $\mu$. The same expansion gives
\[
 (Z-\nu,D(aU-bu))
 \geq a(s^p-s)\,dV+b(1-|Z|)\,d\mu.
\]
Since $s^p-s\geq-\varepsilon/e$ and
$a,b\leq3/2$, this proves \eqref{bvdefect}.

For the Lipschitz case, write $t=|\nabla u|$, so $d\mu=t\,dV$ and
$|Z|=s^\varepsilon$. The two lower bounds above
factor as
\[
 (bs-at)(s^\varepsilon-1)\,dV,\qquad
 (as-bt)(s^\varepsilon-1)\,dV.
\]
For every $K\geq0$ there is $C_K<\infty$ such that
\[
 (s-\lambda)(s^\varepsilon-1)\geq-C_K\varepsilon
 \quad(s\geq0,\ 0\leq\lambda\leq K,\ 0<\varepsilon\leq1).
\]
Indeed, for $0<s\leq1$, the negative part is at most
$s-s^{1+\varepsilon}\leq\varepsilon s|\log s|\leq\varepsilon/e$,
and at $s=0$ the product is nonnegative. For $s>1$, negativity
requires $s<\lambda\leq K$; concavity gives
$s^\varepsilon-1\leq\varepsilon(s-1)$, so the negative part is at
most $\varepsilon K^2$. Apply this bound with $\lambda=at/b$ and
$\lambda=bt/a$, respectively. Since $b/a\leq3$ and
$t\leq\|\nabla u\|_\infty$, both expressions are bounded below by
$-C_{\|\nabla u\|_\infty}\varepsilon\,dV$. This proves
\eqref{lipschitzdefect}.
\end{proof}

\section{Regularity and approximation}\label{bvsection}\label{finitecomparison}

The flux excess in \eqref{bvdefect} is integrated against the full
measure $|Du|$, including the Cantor part. We use a global flux bound
to pass from interior estimates to convergence in $L^1(G,|Du|)$.

\begin{proposition}[Identification of the limit]\label{limitidentification}
Let $G\Subset M$ be a smooth domain with nonempty boundary, and let
$u\in BV(G)\cap L^\infty(G)$ admit a calibration. Suppose that
$U_j\in W^{1,p_j}(G)\cap\Lip_{\mathrm{loc}}(G)$ are $p_j$-IMCFs,
where $p_j\downarrow1$, with continuous fluxes
$Z_j=|\nabla U_j|^{p_j-2}\nabla U_j$. Assume $1\leq u,U_j\leq L$ and
\begin{equation}\label{identificationhypotheses}
 \begin{gathered}
 \sup_j\|Z_j\|_{L^\infty(G)}<\infty,\qquad
 \int_G(|Z_j|-1)_+\,d|Du|\longrightarrow0,\\
 \Tr U_j\longrightarrow\Tr u\quad\text{in }L^1(\partial G).
 \end{gathered}
\end{equation}
Then $U_j\to u$ in $L^{2n/(n-1)}(G)$.
\end{proposition}

\begin{proof}
Discard finitely many terms so that $1<p_j<2$.
Put $\mu=|Du|$, $\varepsilon_j=p_j-1$, $h=\Tr u$, and
$h_j=\Tr U_j$. Set
\[
 A_j=\varepsilon_j|G|+
     \int_G(|Z_j|-1)_+\,d\mu+
     \|h_j-h\|_{L^1(\partial G)}.
\]
The hypotheses give $A_j\to0$. By Lemma~\ref{measuredefect},
\eqref{comparisondefectbounds} holds with
\[
 \eta_j=C\varepsilon_j\,dV+C(|Z_j|-1)_+\,d\mu.
\]
Apply the boundary form of Proposition~\ref{commontruncation} with
$k=0$ and $U=U_j$. Since $h,h_j\geq1$,
\[
 \Tr\bigl(|u-U_j|-\theta(u+U_j)\bigr)_+\leq|h-h_j|.
\]
The uniform flux bound makes the comparison constant independent of
$j$. Thus, with $q=2n/(n-1)$,
\[
 \|u-U_j\|_{L^q(G)}
 \leq2L|G|^{1/q}\theta+C(A_j/\theta)^{1/2}.
\]
Let $j\to\infty$ for fixed $\theta$, and then let $\theta\downarrow0$.
\end{proof}

To prove Theorem~\ref{bvregularity}, it suffices to work on sufficiently small geodesic balls. Set \(\mu=|Du|\) and take such a ball \(G\Subset\Omega\) on which \(r^2=d(x_0,\cdot)^2\) is smooth and has positive definite Hessian.
After a translation, assume $1\leq u\leq L$ on $G$. The trace
$h=\Tr u$ has the same bounds, and $\mu(G)<\infty$.
Smooth boundary data allow us to estimate the flux up to $\partial G$.

For smooth $h_j:\partial G\to[1,L]$ with
\begin{equation}\label{bvtraceapprox}
 h_j\longrightarrow h\quad\text{in }L^1(\partial G),
\end{equation}
choose smooth extensions $H_j:\overline G\to[1,L]$ by extending in a
collar and interpolating with a constant in the interior.

For $\varepsilon=p-1>0$, let $V_{j,p}$ be the positive $p$-harmonic
function in $G$ with boundary value $e^{-h_j/\varepsilon}$, and set
\[
 U_{j,p}=-\varepsilon\log V_{j,p},\qquad
 Z_{j,p}=|\nabla U_{j,p}|^{p-2}\nabla U_{j,p}.
\]
The Dirichlet construction of Section~\ref{dirichletsection} gives
\[
 1\leq U_{j,p}\leq L,\qquad \Tr U_{j,p}=h_j,\qquad
 \diver Z_{j,p}=|\nabla U_{j,p}|^p.
\]
Boundary regularity \cite{Lieberman} gives
$V_{j,p},U_{j,p}\in C^{1,\alpha}(\overline G)$ for fixed $j,p$.
In particular, $Z_{j,p}$ is continuous; the regularity constants may
depend on $j,p$.

For a sequence $p_j\downarrow1$, write $U_j=U_{j,p_j}$ and
$Z_j=Z_{j,p_j}$. Integrals against $\mu$ use the continuous interior
representative of $Z_j$.

The constant $C_j$ in the boundary gradient estimate may grow with
the derivatives of $H_j$. In the flux bound it appears as $C_j^{p-1}$,
so we can keep it bounded by choosing $p_j$ sufficiently close to $1$
after fixing $H_j$.

\begin{lemma}\label{diagonalflux}
There are $p_j=1+\varepsilon_j\downarrow1$ for which
\begin{equation}\label{eq:globalflux}
 \sup_j\norm{Z_j}_{L^\infty(G)}<\infty,
\end{equation}
and
\begin{equation}\label{eq:measuredefect}
 \int_G(|Z_j|-1)_+\,d\mu\longrightarrow0.
\end{equation}
The gradients $\nabla U_j$ are uniformly bounded on each compact
$K\Subset G$.
\end{lemma}

\begin{proof}
Fix $j$ and set $\varepsilon=p-1$. It suffices to show
$\norm{\nabla U_{j,p}}_\infty\leq C_j\varepsilon^{-1}e^{(L-1)/\varepsilon}$,
since the exponential factor becomes $e^{L-1}$ in the flux bound.
The constants $C_j,C'_j,C''_j$ may depend on $(G,g)$, $L$, and $H_j$,
and are independent of $p$ for
$0<\varepsilon\leq1/2$.

Let $\rho=r^2-R^2$ be a defining function for $G$, so that
$\rho<0$ in $G$, $\rho=0$ on $\partial G$, and
$\nabla^2\rho\geq c_g g$ on $\overline G$. Set
$F=e^{-H_j/\varepsilon}$. For $0<\varepsilon\leq1/2$,
\[
 \norm F_{C^2(\overline G)}
 \leq C_j\varepsilon^{-2}e^{-1/\varepsilon}.
\]
Choose $A\leq C'_j\varepsilon^{-2}e^{-1/\varepsilon}$ sufficiently
large that $F+A\rho$ is geodesically convex and $F-A\rho$ is geodesically
concave. A geodesically convex $C^2$ function $Q$ is weakly $p$-subharmonic:
the regularized flux
$(\delta+|\nabla Q|^2)^{(p-2)/2}\nabla Q$ has nonnegative divergence,
and one lets $\delta\downarrow0$.
The analogous concave function is a supersolution. Comparison therefore
gives
\[
 F+A\rho\leq V_{j,p}\leq F-A\rho.
\]
Taking one-sided normal derivatives at the common boundary values, and
using the prescribed tangential derivatives, yields
\[
 \sup_{\partial G}|\nabla V_{j,p}|
 \leq C''_j\varepsilon^{-2}e^{-1/\varepsilon}.
\]
Since $V_{j,p}\geq e^{-L/\varepsilon}$,
\begin{equation}\label{eq:boundarygradient}
 \sup_{\partial G}|\nabla U_{j,p}|
 \leq C''_j\varepsilon^{-1}e^{(L-1)/\varepsilon}.
\end{equation}

Choose $\kappa\geq0$ so that $\operatorname{Ric}\geq-\kappa g$ on $G$.

\begin{claim}
For a $p$-IMCF $U\in C^1(\overline G)$,
\begin{equation}\label{eq:globalgradient}
 \sup_G|\nabla U|
 \leq\max\{\sup_{\partial G}|\nabla U|,\sqrt{(n-1)\kappa}\}.
\end{equation}
\end{claim}

\begin{proof}[Proof of the claim]
Set $f=|\nabla U|^2$. The bound follows if $f$ attains its maximum
on the boundary or if the maximum is zero. At a positive interior
maximum, $U$ is smooth since its gradient is nonzero.
Put $e=\nabla U/|\nabla U|$ and
$A^{ij}=g^{ij}+(p-2)e^ie^j$. The equation is
\[
 A^{ij}U_{ij}=f.
\]
At this maximum, $\nabla^2U(e,\cdot)=0$, hence $\Delta U=f$ and
$|\nabla^2U|^2\geq f^2/(n-1)$. The nonlinear Bochner formula of \cite[Lemma 2.1]{KN}, applied to the linearized operator
\[
\mathcal L\varphi
=\operatorname{div}\bigl(f^{p/2-1}A\nabla\varphi\bigr)
-pf^{p/2-1}\langle\nabla U,\nabla\varphi\rangle,
\]
gives
\[
\mathcal L f
=2f^{p/2-1}
\bigl(|\nabla^2U|^2+
\operatorname{Ric}(\nabla U,\nabla U)\bigr)
+\left(\frac p2-1\right)f^{p/2-2}|\nabla f|^2.
\]
At the interior maximum, \(\nabla f=0\), so the first-order terms vanish. Since \(A\) is positive definite for \(p>1\), we obtain
\[
\begin{aligned}
0
&\ge \frac12 A^{ij}f_{ij} =|\nabla^2U|^2+
\operatorname{Ric}(\nabla U,\nabla U) \ge \frac{f^2}{n-1}-\kappa f.
\end{aligned}
\]
Thus $f\leq(n-1)\kappa$ at the maximum, which gives
\eqref{eq:globalgradient}.
\end{proof}

Apply \eqref{eq:globalgradient} to $U_{j,p}$ and use
\eqref{eq:boundarygradient}. Increasing $C_j\geq1$ to include the
curvature term gives
\begin{equation}\label{eq:fluxpower}
 \norm{Z_{j,p}}_\infty
 =\norm{\nabla U_{j,p}}_\infty^{\varepsilon}
 \leq e^{L-1}(C_j/\varepsilon)^{\varepsilon}.
\end{equation}
Choose $0<\varepsilon_j\leq1/2$ decreasing to zero with
$\varepsilon_j\log C_j\to0$. Evaluating \eqref{eq:fluxpower} at
$p_j=1+\varepsilon_j$ gives \eqref{eq:globalflux}, since
$\varepsilon_j|\log\varepsilon_j|\to0$.

For each compact $K\Subset G$, the interior gradient estimate
\cite[Theorem~1.1]{KN} gives a constant $B_K<\infty$ independent of
$j$ for all sufficiently large $j$. Then $|\nabla U_j|\leq B_K$
on $K$, so
\[
 \sup_K(|Z_j|-1)_+
 \leq \max\{1,B_K\}^{\varepsilon_j}-1\longrightarrow0.
\]
The global flux bound dominates $(|Z_j|-1)_+$ on the finite measure
space $(G,\mu)$. Dominated convergence proves
\eqref{eq:measuredefect} for the full measure $\mu$, including its
singular part.
\end{proof}

\begin{proof}[Proof of Theorem~\ref{bvregularity}]
On each sufficiently small ball $G$ as above,
Lemma~\ref{diagonalflux} and \eqref{bvtraceapprox} give a sequence
satisfying \eqref{identificationhypotheses}.
Proposition~\ref{limitidentification} therefore gives
$U_j\to u$ in $L^{2n/(n-1)}(G)$. The interior gradient bounds give
local equicontinuity. Since the $U_j$ are uniformly bounded,
Arzel\`a--Ascoli on a compact exhaustion of $G$ yields a diagonal
subsequence $U_{j_k}\to v$ locally uniformly. On each relatively
compact geodesically convex ball, the gradient bound gives a common
Lipschitz constant that passes to $v$. Thus $v$ is locally Lipschitz,
and the $L^{2n/(n-1)}$ convergence gives $v=u$ a.e. on $G$.

These representatives agree on overlaps and define a locally
Lipschitz representative on $\Omega$. The calibration reduces to
\eqref{cal}, and the variational inequality following
Definition~\ref{caldef} gives the weak solution property. Identifying
$u$ with the locally Lipschitz limit thus rules out the entire singular
part of $Du$.
\end{proof}

For the Dirichlet problem for least-gradient functions see
\cite{SWZ}, and for their $p$-harmonic approximation see \cite{Juutinen}.

\begin{proposition}\label{sobolevconvergence}
For $G,u,u_p$ as in Theorem~\ref{main}, there is a constant $C$,
depending only on $(G,g)$,
$\operatorname{osc}_G u$, and $\|\nabla u\|_{L^\infty(G)}$, such that
\begin{equation}\label{uniformrate}
 \|u_p-u\|_{C^0(\overline G)}\leq C(p-1)^{1/3}
 \qquad(1<p<2).
\end{equation}
\end{proposition}
\begin{proof}
Translate $u$ and $u_p$ so that $\min_{\overline G}u=1$. Then
$1\leq u,u_p\leq L=1+\operatorname{osc}_G u$. Set $\varepsilon=p-1$
and, for $0<\theta\leq1/2$, put
\[
 w_\theta=|u-u_p|-\theta(u+u_p),\qquad q=\frac{2n}{n-1}.
\]
With $U=u_p$, the Lipschitz bound \eqref{lipschitzdefect} gives
\eqref{comparisondefectbounds} with $\eta=C\varepsilon\,dV$.
For $a=1-\theta$, $b=1+\theta$, and $k\geq0$, both
$au-bu_p-k$ and $au_p-bu-k$ have boundary value $-2\theta u-k<0$.
Continuity makes their positive parts compactly supported, so
Proposition~\ref{commontruncation} gives
\begin{equation}\label{scaledsobolevbound}
 \|(w_\theta-k)_+\|_{L^q(G)}^2
 \leq\frac{C\varepsilon}{\theta}|\{w_\theta>k\}|,
\end{equation}
where $|\cdot|$ denotes Riemannian volume and $C$ depends only on the
data in the proposition.

Set \(\Phi(k)=\int_G(w_\theta-k)_+\,dV\) and \(B=C(\varepsilon/\theta)^{1/2}\), with \(C\) depending on the same data. By the layer-cake formula,
\[
\Phi(k)
=\int_G (w_\theta-k)_+\,dV
=\int_G\int_k^\infty
\mathbf{1}_{\{w_\theta(x)>t\}}\,dt\,dV(x)
=\int_k^\infty |\{w_\theta>t\}|\,dt.
\]
so \(\Phi\) is absolutely continuous and \(-\Phi'(k)=|\{w_\theta>k\}|\) for almost every \(k\). Hölder's inequality and \eqref{scaledsobolevbound} give
\[
\Phi(k)\le B(-\Phi'(k))^{1+1/(2n)}.
\]
Consequently, for almost every $k$ with $\Phi(k)>0$,
\[
\begin{aligned}
-\frac{d}{dk}\Phi(k)^{1/(2n+1)}
&=-\frac{\Phi'(k)}
{(2n+1)\Phi(k)^{2n/(2n+1)}}
\ge \frac{B^{-2n/(2n+1)}}{2n+1}.
\end{aligned}
\]
Since $\Phi(0)\leq B|G|^{1+1/(2n)}$, integration forces $\Phi$ to vanish
by height $(2n+1)B|G|^{1/(2n)}$. Continuity of $w_\theta$ therefore gives
$\|(w_\theta)_+\|_{C^0(\overline G)}\leq C(\varepsilon/\theta)^{1/2}$.
Since $|u_p-u|\leq(w_\theta)_++2L\theta$, we obtain
\[
 \|u_p-u\|_{C^0(\overline G)}
 \leq C\left(\theta+
             \left(\frac{\varepsilon}{\theta}\right)^{1/2}\right).
\]
Taking $\theta=\tfrac12\varepsilon^{1/3}$ gives \eqref{uniformrate}.
The estimate also holds in $L^q(G)$ and $L^1(G)$ since $|G|<\infty$.
It is invariant under common translations of $u,u_p$.
\end{proof}

\begin{proof}[Proof of Theorem~\ref{main}]
Existence is given in Section~\ref{dirichletsection}, and
Proposition~\ref{sobolevconvergence} proves \eqref{mainconvergence}.
For the variation estimate, choose a calibration $\nu$ of $u$ and put
$\varepsilon=p-1$,
$\delta=\|u_p-u\|_\infty$, $s=|\nabla u_p|$, $t=|\nabla u|$, and
$Z_p=|\nabla u_p|^{p-2}\nabla u_p$. The integrating factor gives
$\diver(e^{-u_p}Z_p)=\diver(e^{-u}\nu)=0$.
The common trace makes $\phi=e^{u_p}-e^u$ an admissible
$W^{1,p}_0(G)$ test for both fluxes. Testing the first and applying
H\"older's inequality yields
\[
 A_p:=\int_Gs^p\dd
 =\int_Ge^{u-u_p}Z_p\cdot\nabla u\dd
 \leq e^\delta A_p^{(p-1)/p}
       \left(\int_Gt^p\dd\right)^{1/p},
\]
hence $A_p\leq e^{p\delta}\int_Gt^p\dd$.
Testing the second gives
\[
 \int_Gt\dd=\int_Ge^{u_p-u}\nu\cdot\nabla u_p\dd
 \leq e^\delta\int_Gs\dd.
\]
Since $s\leq s^p+\varepsilon/e$ and
$t^p\leq t+C_{\|\nabla u\|_\infty}\varepsilon$, integration gives
the following bounds with $C$ depending only on $|G|$ and
$\|\nabla u\|_\infty$.
\[
 e^{-\delta}\int_Gt\dd\leq\int_Gs\dd
 \leq e^{p\delta}\left(\int_Gt\dd+C\varepsilon\right)
       +C\varepsilon.
\]
As $\delta\leq\operatorname{osc}_G u$, the difference of total variations
is at most $C(\delta+\varepsilon)$. Proposition~\ref{sobolevconvergence}
proves \eqref{mainstrictbvrate} and strict $BV$ convergence, with the stated
dependence of the constants.
\end{proof}

\begin{proof}[Proof of Corollary~\ref{localapprox}]
By \cite[Lemma~2.3.4]{XuThesis}, the given weak solution admits a
calibration on each smooth $G\Subset\Omega$.
Let $D\Subset\Omega$ lie in a component $\Omega_0$ of
$\Omega$. If $\Omega_0$ is compact and without boundary, weak
minimality with the variation equal to $u-c$ on $\Omega_0$ and
$u$ elsewhere, where $c>0$, gives
$\int_{\Omega_0}|\nabla u|=0$. Thus $u$ is constant on $\Omega_0$,
and one may take $u_p=u$ there.
Otherwise, choose a smooth $G\Subset\Omega_0$ with nonempty boundary
containing $\overline D$ and restrict the Dirichlet approximations.
\end{proof}

\section{Local estimates and boundary values}\label{mainproofs}\label{annularsection}

\begin{corollary}\label{gradientcor}
Suppose $B(x,R)\Subset\Omega$ and $\sec\geq-K$ on $B(x,R)$,
where $K\geq0$. Every locally Lipschitz weak solution satisfies
\begin{equation}\label{gradientbound}
 |\nabla u|\leq C(n,KR^2)R^{-1}
 \quad\text{a.e. on }B(x,R/2).
\end{equation}
\end{corollary}

The sectional curvature hypothesis gives constants uniform as
$p\downarrow1$. The constants in the Ricci-curvature estimate of
\cite{WZ} do not have the uniformity needed for this passage.

\begin{proof}[Proof of Corollary~\ref{gradientcor}]
Let $\Omega_0$ be the connected component of $\Omega$ containing
$x$. If it is compact, the compact-component argument in the proof
of Corollary~\ref{localapprox} shows that $u$ is constant there, so
\eqref{gradientbound} holds. Otherwise, choose a smooth
domain $G$ with nonempty boundary such that
$\overline{B(x,3R/4)}\subset G\Subset B(x,R)$ and use the Dirichlet
approximations from Theorem~\ref{main}. For every $y\in B(x,R/2)$, the ball
$B(y,R/4)$ is contained in $G$. Applied to the positive
$p$-harmonic potentials $v_p$ on these balls, the logarithmic gradient
estimate of \cite[Theorem 1.1]{KN}, with its curvature parameter
$\sqrt K$, gives $|\nabla u_p|\leq C(n,KR^2)R^{-1}$ on
$B(x,R/2)$, uniformly for $p$ sufficiently close to 1. On each
interior compact subset, the gradients are weak-* precompact in
$L^\infty$; the uniform convergence of $u_p$ identifies every such
limit with $\nabla u$ in distributions. The bound therefore passes
to $|\nabla u|$ almost everywhere.
\end{proof}

For weak solutions $u,v\in W^{1,\infty}(G)\cap C(\overline G)$ on a
smooth $G\Subset M$ with nonempty boundary, \cite[Theorem~2.2(i)]{HI}
gives
\[
 \|u-v\|_{C^0(\overline G)}\leq\|u-v\|_{C^0(\partial G)}.
\]
For $u,u_p$ as in Theorem~\ref{main}, let $\widetilde u_p\in
W^{1,p}(G)\cap C(\overline G)$ be a $p$-IMCF with trace error
$\rho_p=\|\widetilde u_p-u\|_{C^0(\partial G)}$. Comparing with
$u_p\pm\rho_p$ gives $\|\widetilde u_p-u_p\|_{C^0(\overline G)}\leq\rho_p$.
By \eqref{mainconvergence}, for $1<p<2$,
\[
 \|\widetilde u_p-u\|_{C^0(\overline G)}
 \leq\rho_p+C(p-1)^{1/3}.
\]

The Dirichlet obstruction in the planar annulus is already noted in
\cite[p.~3072]{XuObstacle}. The next example computes the total
variation of the radial approximations in general dimension.
For the function-obstacle formulation see \cite{MoserObstacle}, and
for $p$-capacitary approximation of the flow with a geometric outer
obstacle see \cite[Theorem~1.11]{BMRSX}.

\begin{example}\label{annulusexample}
Let $G=\{1<r<R\}\subset\mathbb R^n$, where $n\geq2$ and $R>1$.
Put $\alpha=n-1$ and fix $B>\alpha\log R$. For $1<p<n$, let $u_p$ be the
weak solution of
\[
 \begin{cases}
 \Delta_pu_p=|\nabla u_p|^p&\text{in }G,\\
 u_p=0&\text{on }\{r=1\},\\
 u_p=B&\text{on }\{r=R\}.
 \end{cases}
\]
Then, as $p\downarrow1$,
\[
 u_p\longrightarrow v(r)=\alpha\log r
 \quad\text{locally uniformly in }G\text{ and in }L^1(G).
\]
The limit is a calibrated weak solution, but its outer trace satisfies
$v(R)=\alpha\log R<B$. Moreover,
\begin{equation}\label{annulusvariation}
 \lim_{p\downarrow1}\int_G|\nabla u_p|\dd
 =\int_G|\nabla v|\dd
   +|\partial B_R|\,[B-\alpha\log R].
\end{equation}
Thus the convergence is not strict in $BV(G)$, and the prescribed boundary
values are not the trace of a calibrated function in
$W^{1,\infty}(G)\cap C(\overline G)$.
\end{example}

\begin{proof}
Put $\varepsilon=p-1$ and $k_p=(n-p)/\varepsilon$.
The positive radial $p$-harmonic Dirichlet solution and its
logarithmic transform are
\[
 v_p(r)=\frac{r^{-k_p}-R^{-k_p}
       +e^{-B/\varepsilon}(1-r^{-k_p})}{1-R^{-k_p}},
 \qquad u_p=-\varepsilon\log v_p.
\]
For $0<\delta<R-1$, factoring out $r^{-k_p}$ gives
\[
 \frac{v_p(r)}{r^{-k_p}}
 =\frac{1-(r/R)^{k_p}+e^{-B/\varepsilon}(r^{k_p}-1)}
 {1-R^{-k_p}}.
\]
The last term in the numerator tends uniformly to zero on
$[1,R-\delta]$ because $B>\alpha\log R$ and
$\varepsilon k_p=n-p<\alpha$. The remaining terms also converge
uniformly, so $v_p(r)/r^{-k_p}\to1$ there. Hence
\[
 u_p\longrightarrow v(r)=\alpha\log r
 \quad\text{locally uniformly in }G,
 \qquad v(R)=\alpha\log R<B.
\]
The limit $v$ is a weak solution, calibrated by $\partial_r$ since
$\operatorname{div}(\partial_r)=(n-1)/r=|\nabla v|$.
Since $0\leq u_p\leq B$, bounded convergence also gives $u_p\to v$
in $L^1(G)$. In the formula for $v_p$, the coefficient of
$r^{-k_p}$ is positive; hence $v'_p<0$ and $u'_p>0$.
Radial integration by parts gives
\[
 \int_G|\nabla u_p|\dd
 =|\mathbb S^{n-1}|\left(R^{n-1}B
       -\alpha\int_1^R r^{n-2}u_p(r)\,dr\right).
\]
Applying the same identity to $v$ and using bounded convergence yields
\eqref{annulusvariation}.
The excess in \eqref{annulusvariation} is precisely the outer-boundary
area multiplied by the trace discrepancy. It is positive, and hence
the convergence is not strict in $BV(G)$.
If the boundary values were attained by a calibrated
$u\in W^{1,\infty}(G)\cap C(\overline G)$,
Theorem~\ref{main} would give uniform convergence to $u$ on
$\overline G$. The interior limit identifies $u$ with $v$, which
contradicts continuity with boundary value $B$ at $r=R$.
\end{proof}

\ifx\paperaffiliations\empty\else
  \bigskip\noindent\paperaffiliations
\fi
\end{document}